\documentstyle[12pt]{article}
\def\Sym{{\mathrm{Sym}}}
\def\Char{{\mathrm{Char}}}
\def\Shuf{{\mathrm{Shuf}}}
\newtheorem{thm}{Theorem}[subsection]
\newtheorem{lem}[thm]{Lemma}
\newtheorem{defn}[thm]{Definition}
\newtheorem{exm}[thm]{Example}
\begin{document}
\input{amssym}
\renewcommand{\theequation}{\arabic{section}.\arabic{equation}}
\title{Symmetries of 2nd order ODE:
$y''+G(x)y'+H(x)y=0$.}
\author{Mehdi Nadjafikhah\thanks{Department of Mathematics, Iran University of Science and Technology,
Narmak, Tehran, I.R.Iran. e-mail: m\_ nadjafikhah@iust.ac.ir} \and
Seyed-Reza Hejazi\thanks{e-mail: reza\_ hejazi@iust.ac.ir}}
\date{}
\maketitle
\begin{abstract}
This paper is devoted to study the Lie algebra of linear
symmetries of a homogenous 2nd order ODE, by the method of
Kushner, Lychagin and Robstov \cite{[1]}.
\end{abstract}
\noindent {\bf A.M.S. 2000 Subject Classification:} {\em 11Dxx,
32Wxx , 76Mxx.}
\section*{Introduction}
Symmetries of differential equations make a magnificent portion in
theory of differential equations, and there are so much researches
in this object. Here we are going to decompose the structure of
Lie algebra of linear symmetries of ODE, $y''+G(x)y'+H(x)y=0,$
where $G$ and $H$ are smooth functions of $x$, to two subalgebras
which are called  \textbf{even} and \textbf{odd} symmetries of
ODE. First the meaning of even and odd symmetries for a
differential operator is given, next we will find these two
concepts for the differential operator corresponding to ODE
instead of the equation itself. Reader is referred to \cite{[2]}
and \cite{[3]} for some fundamental contexts in geometry of
manifolds and their applications to theory of differential
equations.
\section{Symmetries of ODEs}
Consider a general $n-$th order differential equation $\Delta=0$
which is defined on $n-$th jet space of $p$ independent and $q$
dependent variables. As we know a  \textit{symmetry} of the system
of above differential equation means a point (or contact)
transformation which maps solutions to solutions. In the case of
point transformation, the infinitesimal generator {\bf v} from a
Lie algebra ${\goth g}$ corresponding to group transformation
makes a symmetry of $\Delta=0$, if its $n-$th prolongation
annihilate $\Delta$, i.e., ${{\bf v}}^{(n)}(\Delta)=0.$ See
\cite{[2]} and \cite{[3]} for more details of symmetries of
differential equations. It is noteworthy that all manifolds,
vector fields, differential forms and... are seem to be smooth in
the sequel.
\subsection{Generating Functions}
Let us consider an ODE of $(k+1)-$th order which is resolved with
respect to the highest derivative:
$y^{(k+1)}=F(x,y,y',...,y^{(k)})$. This equation determines a
one-dimensional distribution on the $k-$th jet space with one
independent variable $x$ with coordinate $(x,y=p_0,p_1,...,p_k)$,
which is generated by the vector field $${\cal
D}=\frac{\partial}{\partial x}+p_1\frac{\partial}{\partial
p_0}+\cdots+p_k\frac{\partial}{\partial
p_{k-1}}+F\frac{\partial}{\partial p_k},$$ or by the contact
differential 1-forms
$$\omega^1=dp_0-p_1dx,\cdots,\hspace{0.5cm} \omega^k=dp_{k-1}-p_kdx,\hspace{0.5cm} \omega^{k+1}=dp_k-Fdx.$$
Consider a vector field $X$ on manifold $M$, $X$ is called a
\textit{symmetry} of the distribution $P$, if the distribution is
invariant under the flow of $X$. Denote by $\Sym(P)$ the set of
all symmetries id $P$. If $X$ belongs to $P$ then it is called a
\textit{characteristic symmetry} and the set of all characteristic
symmetries of $X$ is denoted by $\Char(P)$ which makes an ideal of
$\Sym(P)$.
\begin{defn}
The quotient Lie algebra
\begin{eqnarray*}
\Shuf(P)=\Sym(P)/\Char(P),
\end{eqnarray*}
is called the set of \textit{shuffling symmetries} of $P$.
\end{defn}
Therefore any shuffling symmetry $S\in\Shuf(P)$ has a unique
representative of the form
$$S=f\frac{\partial}{\partial p_0}+{\mathcal{D}}(f)\frac{\partial}{\partial p_1}+\cdots+{\mathcal{D}}^k(f)\frac{\partial}{\partial p_k},$$
where $f$ is a smooth functions of $(x,p_0,p_1,...,p_k)$ and
${\mathcal{D}}^i={\mathcal{D}}({\mathcal{D}}^{i-1})$, for the
reason see \cite{[1]}. The function $f$ is called a
\textit{generating function} of the symmetry $S$ and we write
$S_f$ instead of $S$. Therefore, $S_f$ is a shuffling symmetry of
the ODE if and only if the generation function $f$ satisfies the
following \textit{Lie equation}:
$${\mathcal D}^{k+1}(f)-\sum_{i=0}^{k}\frac{\partial F}{\partial p_i}{\mathcal D}^i(f)=0.$$
\par
Let us denote by $\Delta_F:C^{\infty}(\Bbb{R}^{k+2})\rightarrow
C^{\infty}(\Bbb{R}^{k+2})$, the following linear $k-$th order
scalar differential operator:
$$\Delta_F={\mathcal D}^{k+1}-\sum_{i=0}^{k}\frac{\partial F}{\partial p_i}{\mathcal D}^i,$$
which is called the linearization of
$y^{(k+1)}=F(x,y,y',...,y^{(k)})$.
\begin{thm}
There exist the isomorphism $\Shuf(P)\cong\ker\Delta_F$ between
solutions of the Lie equation and shuffling symmetries.
\end{thm}
\par
The $\Shuf(P)$ is a Lie algebra for any distribution $P$ with
respect to the \textit{Poisson-Lie bracket}, which is defined in
the following way:
\begin{eqnarray*}
[S_f,S_g]:=S_{[f,g]}=\sum_{i=0}^{k}\bigg({\mathcal
D}^i(f)\frac{\partial g}{\partial p_i}-{\mathcal
D}^i(g)\frac{\partial f}{\partial p_i}\bigg)
\end{eqnarray*}
for any $f,g\in\ker\Delta_F.$

\begin{exm}
Functions $f=a(x,p_0)p_1+b(x,p_0)$ are generating functions of the
vector fields on $\Bbb{R}^2$ of the form
$b(x,p_0)\frac{\partial}{\partial
p_0}-a(x,p_0)\frac{\partial}{\partial x}.$
\end{exm}
\subsection{Linear Symmetries} A shuffling symmetry $S_f$ is called
a \textit{linear symmetry}, if the generating function $f$ is
linear in $p_0,...,p_k$, i.e., $f=b_0(x)p_0+\cdots+b_k(x)p_k.$
With any linear symmetry we associate a linear operator
$\Delta_f=b_0+\cdots+b_k\partial^k$, where $\partial=d/dx$, and we
rewrite the Lie equation for linear symmetries in terms of the
algebra of linear differential operators.
\begin{lem}
For any linear differential operator $A=a_0+\cdots+a_n\partial^n$
and $L=l_0+\cdots+l_k\partial^k+\partial^{k+1}$ there are unique
differential operators $C_A$ and $R_A$ of order $\leq n-k-1$ and
$\leq k$ respectively such that $A=C_A\circ L+R_A.$
\end{lem}
Here we have a very important theorem:
\begin{thm}\cite{[1]}
A differential operator $\Delta_f=b_0+\cdots+b_k\partial^k$
corresponds to a shuffling symmetry $f=b_0p_0+\cdots+b_kp_k$ of
the linear differential equation $L(h)=0$, where
$L=A_0+\cdots+A_k\partial^k+\partial^{k+1}$, if and only if there
is a differential operator $\nabla_f$ of order k and such that
$L\circ\Delta_f=\nabla_f\circ L$. Moreover, the commutator $[f,g]$
of linear symmetries corresponds to the remainder R of division
$[\Delta_f,\Delta_g]$ by L; that is,
$R_{[\Delta_f,\Delta_g]}=\Delta_{[f,g]}.$
\end{thm}
\par
Denote by ${\frak B}(L)$ the Lie algebra of all differential
operators $\Delta$ such that $L\circ\Delta=\nabla\circ L,$ for
some uniquely determined differential operator $\nabla.$ If
$\Sym(L)$ denote the Lie algebra of linear symmetries of
differential operator $L$, then we have
\begin{thm}
\begin{itemize}
\item[1.] If $\Delta\in {\frak B}(L)$ then $R_{\Delta}\in\Sym(L).$
\item[2.] The residue map $R:{\frak B}(L)\rightarrow \Sym(L),$ is a
Lie algebra homomorphism.
\end{itemize}
\end{thm}
\section{Linear Symmetries of Operators} The differential
operator
$$L^T=(-1)^{k+1}\partial^{k+1}+\sum_{i=0}^k(-1)^i\partial^i\circ A_i,$$
is said to be \textit{adjoint} to the operator
\begin{eqnarray}
L=\partial^{k+1}+\sum_{i=0}^kA_i\partial^i.\label{eq:1}
\end{eqnarray}
\par
A differential operator $L$ is said to be \textit{self-adjoint} if
$L^T=L$ and \textit{skew-adjoint} if $L^T=-T$.\\
\par
The correspondence $\Delta_{f}\leftrightarrow\nabla_{f}^T$
establishes an isomorphism between linear symmetries of the
differential equation $L(h)=0$ and linear symmetries of the
adjoint equation $L^T(h)=0.$
\subsection{$\Bbb{Z}_2-$ Grading on ${\frak B}(L)$}
Let us now assume that $L$ is self-adjoint or skew-adjoint. Then
if $\Delta\in{\frak B}(L)$ so $\nabla^T$ does.\\
using the involution we can decompose ${\frak B}(L)$ is to the
direct some: $${\frak B}(L)={\frak B}_0(L)\oplus{\frak B}_1(L)$$
where
\begin{eqnarray*}
{\frak B}_0(L)&=&\Big\{\Delta:L\circ\Delta=-\Delta^T\circ L\Big\},\\
{\frak B}_1(L)&=&\Big\{\Delta:L\circ\Delta=\Delta^T\circ L\Big\}.
\end{eqnarray*}
We will define $\Bbb{Z}_2-$parity
$\varepsilon(\Delta)=0\in\Bbb{Z}_2$ for $\Delta\in{\frak B}_0(L)$
and $\varepsilon(\Delta)=1\in\Bbb{Z}_2$ for $\Delta\in{\frak
B}_1(L)$, and will consider the above decomposition as
$\Bbb{Z}_2-$grading on ${\frak B}(L)$.
\begin{thm}
Let $L$ be a self or skew-adjoint differential operator.
\begin{itemize}
\item[1.] Then the commutator of operators determines a Lie algebra
structure on ${\frak B}(L)$, such that $$[\Delta_a,\Delta_b]\in
{\frak B}_{a+b}(L)$$ if $\Delta_a\in{\frak B}_a(L)$,
$\Delta_b\in{\frak B}_b(L)$, $a$, $b\in\Bbb{Z}_2.$
\item[2.] Let $\Sym(L)$ be the Lie algebra of Linear symmetries of
operator $L$, and $\Sym_a(L)=R({\frak B}_a(L))$ for
$a\in\Bbb{Z}_2$ and $\Sym_b(L)=R({\frak B}_b(L))$ for
$b\in\Bbb{Z}_2$. Then $\Sym(L)=\Sym_0(L)\oplus\Sym_1(L)$ and
$$[\Sym_a(L),\Sym_b(L)]\subset\Sym_{a+b}(L).$$
\end{itemize}
\end{thm}
\par
We call elements of $\Sym_0(L)$ by \textit{even symmetries} and
elements of $\Sym_1(L)$ by \textit{odd symmetries} of the equation
$L(h)=0.$
\subsection{Symmetries of Operator $\partial^2+G(x)\partial+H(x)$ }
In this part we apply the results to the corresponding operator of
2nd ODE
\begin{eqnarray}
y''+G(x)y'+H(x),\label{eq:3}
\end{eqnarray}
 It is easy to see that in order two one has only self
adjoint operator, thus the following operator is self adjoint and
consequently, we will work on
\begin{eqnarray}
\partial^2+G(x)\partial+H(x),\label{eq:4}
\end{eqnarray}
instead of equation {\ref{eq:3}}. This operator is self adjoint.
Therefore the algebra of linear symmetries in $\Bbb{Z}_2-$graded.
\par
Let us begin with $\Sym _0(L)$. If $\Delta=A_0+A_1\partial\in\Sym
_0(L)$ then we have $L\circ\Delta=-\Delta^T\circ L$. If
$\Delta^T=A_0-A_1'-A_1\partial$ then
\begin{eqnarray*}
L\circ\Delta   &=&A_1\partial^3+\big(A_0+2A_1'+A_1G\big)\partial^2\\
               &&+\Big[2A_0'+A_1''+G\big(A_0+A_1'\big)+A_1H\Big]\partial\\
               & &+A_0''+GA_0'+HA_0,\\
\Delta^T\circ L
               &=&-A_1\partial^3+\big(A_0-A_1'-A_1G\big)\partial^2 \\
               &&+\Big[G\big(A_0-A_1'\big)-A_1\big(G+G'\big)\Big]\partial\\
               & &+H\big(A_0-A_1'\big)-A_1H'.
\end{eqnarray*}
Therefore, $\Delta\in \Sym _0(L)$ implies $A_0=-\frac{1}{2}A_1'$
and the function $A_1=w$ should satisfy the following differential
equation:
\begin{eqnarray}
w'''+\Big(2H-G^2-2G'\Big)w'+\Big(H'-GG'-G''\Big)w=0.\label{eq:2}
\end{eqnarray}
We denote the differential operator corresponding to (\ref{eq:2})
by:
\begin{eqnarray}
\widetilde{L}=\partial^3+\Big(2H-G^2-2G'\Big)\partial+\Big(H'-GG'-G''\Big).\label{eq:5}
\end{eqnarray}
If $\Delta\in \Sym _1(L)$ then $L\circ\Delta=\Delta^T\circ L$ and
we obtain $A_1=0$ and $A_0'=0$. Thus $\Delta\in \Sym _1(L)$ is and
only if $\Delta$ proportional to the identity operator. Finally we
have the following theorem which is a generalization of the
theorem 2.5.1 of \cite{[1]}.
\begin{thm}
The Lie algebra of linear symmetries of the differential operator
(\ref{eq:5}) has the following description:
\begin{itemize}
\item[1.] $\Sym_0(L)=\Big\{-\frac{1}{2}w'+w\partial:
\widetilde{L}(w)=0\Big\}.$
\item[2.] $\Sym _1(L)=\Bbb{R}.$
\end{itemize}
\end{thm}
\begin{exm}
Let us suppose that in (\ref{eq:4}) we have $G(x)=0$. Thus the new
operator is called \textbf{Schr\"{o}dinger operator}. It is
possible to see that the even symmetries of Schr\"{o}dinger
operator is isomorphic to the Lie algebra ${\goth sl}(2)$. And of
course the Schr\"{o}dinger operator does not have any nontrivial
odd symmetries.
\end{exm}
\end{document}